\documentclass{amsart}
\usepackage{amssymb,amsthm,amsmath,epsfig,xypic}
\usepackage{calc}
\usepackage{latexsym}
\usepackage{amscd,amssymb,subfigure,hyperref}
\usepackage[arrow,matrix,graph,frame,poly,arc,tips]{xy}
\begin{document}

\newcommand{\mmbox}[1]{\mbox{${#1}$}}
\newcommand{\affine}[1]{\mmbox{{\mathbb A}^{#1}}}
\newcommand{\Ann}[1]{\mmbox{{\rm Ann}({#1})}}
\newcommand{\caps}[3]{\mmbox{{#1}_{#2} \cap \ldots \cap {#1}_{#3}}}
\newcommand{\N}{{\mathbb N}}
\newcommand{\Z}{{\mathbb Z}}
\newcommand{\R}{{\mathbb R}}
\newcommand{\A}{{\mathcal A}}
\newcommand{\B}{{\mathcal B}}
\newcommand{\C}{{\mathbb K}}
\newcommand{\PP}{{\mathbb P}}
\newcommand{\cO}{{\mathcal O}}
\newcommand{\Tor}{\mathop{\rm Tor}\nolimits}
\newcommand{\ot}{\mathop{\rm OT}\nolimits}
\newcommand{\ao}{\mathop{\rm AOT}\nolimits}
\newcommand{\Ext}{\mathop{\rm Ext}\nolimits}
\newcommand{\Hom}{\mathop{\rm Hom}\nolimits}
\newcommand{\Cl}{\mathop{\rm Cl}\nolimits}
\newcommand{\im}{\mathop{\rm Im}\nolimits}
\newcommand{\rank}{\mathop{\rm rank}\nolimits}
\newcommand{\codim}{\mathop{\rm codim}\nolimits}
\newcommand{\supp}{\mathop{\rm supp}\nolimits}
\newcommand{\CB}{Cayley-Bacharach}
\newcommand{\HF}{\mathrm{HF}}
\newcommand{\HP}{\mathrm{HP}}
\newcommand{\Pic}{\mathrm{Pic}}
\newcommand{\coker}{\mathop{\rm coker}\nolimits}
\sloppy
\newtheorem{defn0}{Definition}[section]
\newtheorem{prop0}[defn0]{Proposition}
\newtheorem{conj0}[defn0]{Conjecture}
\newtheorem{thm0}[defn0]{Theorem}
\newtheorem{lem0}[defn0]{Lemma}
\newtheorem{corollary0}[defn0]{Corollary}
\newtheorem{example0}[defn0]{Example}
\newcommand\res[1]{{\lower1pt\hbox{$|$}}_{\raise.5pt\hbox{${\scriptstyle #1}$}}}

\newenvironment{defn}{\begin{defn0}}{\end{defn0}}
\newenvironment{prop}{\begin{prop0}}{\end{prop0}}
\newenvironment{conj}{\begin{conj0}}{\end{conj0}}
\newenvironment{thm}{\begin{thm0}}{\end{thm0}}
\newenvironment{lem}{\begin{lem0}}{\end{lem0}}
\newenvironment{cor}{\begin{corollary0}}{\end{corollary0}}
\newenvironment{exm}{\begin{example0}\rm}{\end{example0}}

\newcommand{\msp}{\renewcommand{\arraystretch}{.5}}
\newcommand{\rsp}{\renewcommand{\arraystretch}{1}}

\newenvironment{lmatrix}{\renewcommand{\arraystretch}{.5}\small
 \begin{pmatrix}} {\end{pmatrix}\renewcommand{\arraystretch}{1}}
\newenvironment{llmatrix}{\renewcommand{\arraystretch}{.5}\scriptsize
 \begin{pmatrix}} {\end{pmatrix}\renewcommand{\arraystretch}{1}}
\newenvironment{larray}{\renewcommand{\arraystretch}{.5}\begin{array}}
 {\end{array}\renewcommand{\arraystretch}{1}}

\def \a{{\mathrel{\smash-}}{\mathrel{\mkern-8mu}}
{\mathrel{\smash-}}{\mathrel{\mkern-8mu}} {\mathrel{\smash-}}{\mathrel{\mkern-8mu}}}

\title[Toric Hirzebruch-Riemann-Roch via Ishida's theorem on the Todd genus]%
{Toric Hirzebruch-Riemann-Roch via \\Ishida's theorem on the Todd genus}
\author{Hal Schenck}
\thanks{Schenck supported by NSF 1068754, NSA H98230-11-1-0170}
\address{Schenck: Mathematics Department \\ University of
 Illinois \\
   Urbana \\ IL 61801\\USA}
\email{schenck@math.uiuc.edu}

\subjclass[2000]{14M25, 14C40} \keywords{Toric variety, Chow ring, cohomology}

\begin{abstract}
\noindent We give a simple proof of the Hirzebruch-Riemann-Roch theorem
for smooth complete toric varieties, based on Ishida's result in \cite{I}
that the Todd genus of a smooth complete toric variety is one.
\end{abstract}
\maketitle
\vskip -1.5in
\section{Introduction}\label{sec:one}
The Hirzebruch-Riemann-Roch theorem relates the Euler characteristic
of a coherent sheaf ${\mathcal F}$ on a smooth complete $n-$dimensional variety 
$X$ to intersection theory, via the formula
\begin{equation}
\label{HRReqn}
\chi(\mathcal F) = \int ch({\mathcal F})Td({\mathcal T}_X).
\end{equation}
In \cite{BV}, Brion-Vergne prove an equivariant Hirzebruch-Riemann-Roch
theorem for complete simplicial toric varieties. If the toric variety 
is actually smooth, it is possible to derive \eqref{HRReqn} from their
result. In this note, we give a simple direct proof of (\ref{HRReqn}) when
$X$ is a smooth complete toric variety. Such a variety is determined
by a smooth complete rational polyhedral fan 
$\Sigma \subseteq N_{\mathbb R}$, where 
$N \simeq \Z^n$ is a lattice; we write $X$ for the associated toric
variety $X_\Sigma$. We will make use of the following standard facts
about toric varieties. First, 
\begin{equation}
Td(X_\Sigma) = \prod\limits_{\rho \in \Sigma(1)} \frac{D_\rho}{1-e^{-D_\rho}},
\end{equation}
where $\Sigma(k)$ denotes the set of $k$-dimensional 
faces of $\Sigma$. For $\tau \in \Sigma(k)$ there is an 
associated torus invariant orbit $O(\tau)$, and we use 
$V(\tau)$ to denote the orbit
closure $\overline{O(\tau)}$, which has dimension $n-k$. 
A key fact is that (see \cite{CLS}, Proposition~3.2.7) 
\begin{equation}
\label{stareqn}
V(\tau) = \overline{O(\tau)} \simeq X_{\mathrm{Star}(\tau)}.
\end{equation}
Since $\Sigma$ is smooth, all orbits are also smooth, and
if $\rho_i, \rho_j$ are distinct elements of $\Sigma(1)$, then 
(see \cite{CLS}, Lemma~12.5.7)
\[
[D_{\rho_i}\res{V(\rho_j)}] = \begin{cases}
                       V(\tau) & \tau = \rho_i + \rho_j \in \Sigma\\
                        0 & \rho_i , \rho_j \text{ are not both in any cone in } \Sigma.
                       \end{cases}
\]
The final ingredient we need is a result of Ishida: building 
on work of Brion \cite{B}, in \cite{I} Ishida shows 
that (\ref{HRReqn}) holds for the structure sheaf of a 
smooth complete toric variety $X$: 
\begin{equation}
\label{chiOX}
1 = \int Td({\mathcal T}_X) = \Big[\prod\limits_{\rho \in \Sigma(1)} \frac{D_\rho}{1-e^{-D_{\rho}}}\Big]_n.
\end{equation}

\section{The proof}\label{sec:two}
For a smooth complete toric variety, any coherent sheaf has 
a resolution by line bundles \cite{cox}, so it suffices to 
consider the case ${\mathcal F} = {\mathcal O}_X(D)$. 
Let $X = X_\Sigma$, and recall that $\Pic(X)$ is generated by 
the classes of the divisors $D_{\rho}$, $\rho \in \Sigma(1)$.
We will show that if (\ref{HRReqn}) holds for a divisor $D$, then 
it also holds for $D + D_{\rho}$ and $D-D_{\rho}$, for any 
$\rho \in \Sigma(1)$. We begin with the case $D-D_{\rho}$, 
and induct on the dimension of $X$.

A smooth complete toric variety of dimension one is 
simply ${\mathbb P}^1$, so the base case holds by Riemann-Roch
for curves. Suppose the theorem holds for all
smooth complete fans of dimension $< n$, and let $\Sigma$ be a smooth complete
fan of dimension $n$. When $D = 0$ the result holds by Ishida's theorem. 
Let $\rho \in \Sigma(1)$, and partition the rays of $\Sigma$ as 
\[
\Sigma(1) = \rho \cup \Sigma'(1) \cup \Sigma''(1),
\]
where the rays in $\Sigma'(1)$ are in one to one correspondence with the rays of
the fan $\mathrm{Star}(\rho)$. Let $X' = X_{\mathrm{Star}(\rho)} \simeq V(\rho)$. 
Tensoring the standard exact sequence
\[
0 \longrightarrow \cO_{X}(-D_{\rho})  \longrightarrow \cO_{X}
\longrightarrow \cO_{X'} \longrightarrow 0
\]
with $\cO_{X}(D)$ yields the sequence
\[
0 \longrightarrow \cO_{X}(D-D_{\rho})  \longrightarrow \cO_{X}(D)
\longrightarrow \cO_{X'}(D) \longrightarrow 0.
\]
From the additivity of the Euler characteristic, we have
\[
\chi(\cO_{X}(D)) - \chi(\cO_{X}(D-D_{\rho})) = \chi(\cO_{X'}(D)).
\]
Our hypotheses imply that
\begin{align*}
\int_{X'} e^D Td({\mathcal T}_{X'}) & = \chi(\cO_{X'}(D))\\
\int_X e^D Td({\mathcal T}_{X})& = \chi(\cO_{X}(D)),
\end{align*}
so it suffices to show that
\begin{equation}
\label{stepOne}
\begin{aligned}
\int_{X'} ch(D) Td({\mathcal T}_{X'}) & = \int_X (e^D -e^{D-D_{\rho}}) Td({\mathcal T}_{X})\\
                                        &=  \int_X e^D \Big( \frac{1-e^{-D_{\rho}}}{D_{\rho}} \Big) D_{\rho} Td({\mathcal T}_{X})
\end{aligned}
\end{equation}
Break the Todd class of $X$ into two parts:
\[
Td({\mathcal T}_{X}) = \prod\limits_{\gamma \in \Sigma'(1) \cup \rho }\frac{D_{\gamma}}{1-e^{-D_{\gamma}}} \cdot \prod\limits_{\gamma \in \Sigma''(1) }\frac{D_{\gamma}}{1-e^{-D_{\gamma}}}
\]
In (\ref{stepOne}), the term $\frac{1-e^{-D_{\rho}}}{D_{\rho}}$ cancels 
with the corresponding term in $Td({\mathcal T}_{X})$, 
so that 
\begin{equation}
\label{stepTwo}
\begin{aligned}
\int_X e^D \Big( \frac{1-e^{-D_{\rho}}}{D_{\rho}}\Big) D_{\rho} Td({\mathcal T}_{X}) &=
\int_X e^D D_{\rho} \prod\limits_{\gamma \in \Sigma'(1) \cup \Sigma''(1)}\frac{D_{\gamma}}{1-e^{-D_{\gamma}}} \\
&= \int_X e^D D_{\rho} \prod\limits_{\gamma \in \Sigma'(1)}\frac{D_{\gamma}}{1-e^{-D_{\gamma}}}. 
\end{aligned}
\end{equation}
The second equality follows since $D_{\rho} \cdot D_{\gamma}=0$ if 
$\gamma \in \Sigma''(1)$. By smoothness, all intersections are either zero or one,
and thus 
\begin{align*}
\int_X e^D D_{\rho} \prod\limits_{\gamma \in \Sigma'(1)}\frac{D_{\gamma}}{1-e^{-D_{\gamma}}} &=\Big[ e^D D_{\rho} \prod\limits_{\gamma \in \Sigma'(1)}\frac{D_{\gamma}}{1-e^{-D_{\gamma}}} \Big]_n \\
 &=\Big[e^{D\res{V(\rho)}} \prod\limits_{\gamma \in \Sigma'(1)}\frac{D_{\gamma}}{1-e^{-D_{\gamma}}} \Big]
_{n-1}\\
                 &= \int_{X'}e^{D} \cdot  Td({\mathcal T}_{X'}).
\end{align*}
This proves the result for $D-D_{\rho}$. For $D+D_{\rho}$, the result follows
using the substitution $e^{D_{\rho}}-1 = e^{D_{\rho}}(1-e^{-D_{\rho}}).$
\vskip .1in\noindent{\bf Question} Ishida's proof (\ref{chiOX}) is not easy; 
does there exist a simple proof of (\ref{chiOX})?
\vskip .1in
\noindent{\bf Acknowledgements} I thank David Cox for pointing out
Ishida's result to me, and the referee for useful comments.
\bibliographystyle{amsalpha}

\end{document}